\documentclass{article}

\newtheorem{theorem}{Theorem}[section]
\newtheorem{corollary}[theorem]{Corollary}

\newtheorem{lemma}[theorem]{Lemma}
\newtheorem{proposition}[theorem]{Proposition}

\newtheorem{remark}{Remark}

\evensidemargin\oddsidemargin
\begin{document}

\author{Vadim E. Levit and Eugen Mandrescu \\
Department of Computer Science\\
Holon Academic Institute of Technology\\
52 Golomb Str., P.O. Box 305\\
Holon 58102, ISRAEL\\
\{levitv, eugen\_m\}@barley.cteh.ac.il}
\date{}
\title{Combinatorial Properties of the Family of Maximum Stable Sets of a Graph }
\maketitle

\begin{abstract}
The \textit{stability number }$\alpha (G)$ of a graph $G$ is the size of a 
\textit{maximum stable set} of $G$, $core(G)=\cap \{S:S\ $\textit{is a
maximum stable in} $G\}$, and $\xi (G)=\left| core(G)\right| $. In this
paper we prove that for a graph $G$ without isolated vertices, the following
assertions are true: ($\mathit{i}$) if $\xi (G)\leq 1$, then $G$ is
quasi-regularizable; ($\mathit{ii}$) if $G$ is of order $n$ and $\alpha
(G)>(n+k-1)/2$, for some $k\geq 1$, then $\xi (G)\geq k+1$, and $\xi (G)\geq
k+2$, whenever $n+k-1$ is even. The last finding is a strengthening of a
result of Hammer, Hansen, and Simeone, which states that $\alpha (G)>n/2$
implies $\xi (G)\geq 1$. In the case of K\"{o}nig-Egerv\'{a}ry graphs, i.e.,
for graphs enjoying $\alpha (G)+\mu (G)=n$, where $\mu (G)$ is the maximum
size of a matching of $G$, we prove that $\left| core(G)\right| >\left|
N(core(G))\right| $ is a necessary and sufficient condition for $\alpha
(G)>n/2$. Moreover, for bipartite graphs without isolated vertices, $\xi
(G)\geq 2$ is equivalent to $\alpha (G)>n/2$. We also show that Hall's
marriage Theorem is valid for K\"{o}nig-Egerv\'{a}ry graphs, and, it is
sufficient to check Hall's condition only for one specific stable set,
namely, for $core(G)$.
\end{abstract}

\section{Introduction}

Throughout this paper $G=(V,E)$ is a simple (i.e., a finite, undirected,
loopless and without multiple edges) graph with vertex set $V=V(G)$ and edge
set $E=E(G).$ If $X\subset V$, then $G[X]$ is the subgraph of $G$ spanned by 
$X$. By $G-W$ we mean either the subgraph $G[V-W]$ , if $W\subset V(G)$, or
the partial subgraph $H=(V,E-W)$ of $G$, for $W\subset E(G)$. Anyway, we use 
$G-w$, whenever $W$ $=\{w\}$. If $A,B$ $\subset V$ and $A\cap B=\emptyset $,
then $(A,B)$ stands for the set $\{e=ab:a\in A,b\in B,e\in E\}$. A stable
set $S$ of maximum size will be referred as to a \textit{maximum stable set}
of $G$, and $\alpha (G)=\left| S\right| $ is the \textit{stability number }%
of $G$. Let $\Omega (G)$ and $core(G)$ denote respectively the sets $\{S:S$ 
\textit{is a maximum stable set of} $G\}$ and $\cap \{S:S\in \Omega (G)\}$,
while $\xi (G)=\left| core(G)\right| $. Clearly, any isolated vertex of a
graph $G$ is contained in $core(G)$. Let us define $isol(G)$ as the set of
isolated vertices of $G$. The neighborhood of a vertex $v\in V$ is the set $%
N(v)=\{w:w\in V$ \textit{and} $vw\in E\}$, and for $A\subset V,N(A)=\cup
\{N(v):v\in A\}$, while $N[A]=A\cup N(A)$. By $C_{n}$, $K_{n}$, $P_{n}$ we
denote the chordless cycle on $n\geq $ $4$ vertices, the complete graph on $%
n\geq 1$ vertices, and respectively the chordless path on $n\geq 3$ vertices.

A \textit{matching} is a set of non-incident edges of $G$; a matching of
maximum cardinality $\mu (G)$ is a \textit{maximum matching}, and a \textit{%
perfect matching} is a matching covering all the vertices of $G$. $G$ is a 
\textit{K\"{o}nig-Egerv\'{a}ry graph }provided $\alpha (G)+\mu (G)=\left|
V(G)\right| $, \cite{dem}, \cite{ster}. According to a well-known result of
K\"{o}nig, \cite{koen}, and Egerv\'{a}ry, \cite{eger}, any bipartite graph
enjoys this property.

A graph $G$ is $\alpha ^{+}$-\textit{stable} if $\alpha (G+e)=\nolinebreak
\alpha (G)$, for any edge $e\in E(\overline{G})$, where $\overline{G}$ is
the complement of $G$, \cite{gun}. The following characterization of $\alpha
^{+}$-stable graphs, without any referring to this notion, had been proved
in \cite{hayn} three years before the above definition was proposed.

\begin{theorem}
\label{th1}\cite{hayn} A graph $G$ is $\alpha ^{+}$-stable if and only if $%
\xi (G)\leq 1$.
\end{theorem}

\begin{figure}[h]
\setlength{\unitlength}{1.0cm} 
\begin{picture}(5,1.6)\thicklines

  \multiput(4,0.5)(1,0){3}{\circle*{0.29}}

  \put(5,1.5){\circle*{0.29}}
  \put(4,0.5){\line(1,0){2}}
  \put(5,0.5){\line(0,1){1}}
  \put(4,0.5){\line(1,1){1}}
  \put(6,0.5){\line(-1,1){1}}

  \put(5,0){\makebox(0,0){$(a)$}} 

  \multiput(8,0.5)(1,0){3}{\circle*{0.29}}

  \put(8,0.5){\line(1,0){2}}
  \put(9,1.5){\circle*{0.29}}
  \put(9,0.5){\line(0,1){1}}
  \put(9,1.5){\line(1,-1){1}}

  \put(9,0){\makebox(0,0){$(b)$}} 
 
 \end{picture}
\caption{{Two non-$\alpha _{0}^{+}$-stable graphs: (\textit{a}) $K_{4}-e$; (%
\textit{b}) $K_{3}+e$.}}
\label{fig1}
\end{figure}
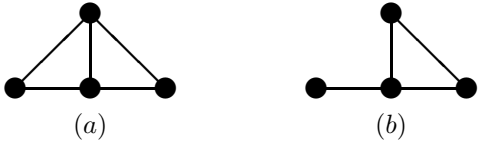

This result motivates that a graph $G$ is referred to as: ($\mathit{i}$) $%
\alpha _{0}^{+}$-stable, if $\xi (G)=0$, and ($\mathit{ii}$) $\alpha
_{1}^{+} $-stable provided $\xi (G)=1$, \cite{levm3}. For instance, $C_{4}$
is $\alpha _{0}^{+}$-stable, $K_{3}+e$ is $\alpha _{1}^{+}$-stable, and the
diamond, i.e., the graph $K_{4}-e$, is not $\alpha ^{+}$-stable (see Figure 
\ref{fig1}).

A graph $G$ is \textit{quasi-regularizable} if one can replace each edge of $%
G$ with a non-negative integer number of parallel copies, so as to obtain a
regular multigraph of degree $\neq 0$, \cite{berge2}. For instance, the
diamond is quasi-regularizable, while $P_{3}$ is not quasi-regularizable.
Clearly, any quasi-regularizable graph has no isolated vertices. Moreover, a
disconnected graph is quasi-regularizable if and only if any of its
connected components is a quasi-regularizable graph.

In this paper we analyze the relationship between $\alpha (G)$ and $\xi (G)$%
. We show that if $G$ has $\left| isol(G)\right| \neq 1$ and $\alpha
(G)>(\left| V(G)\right| +k-1)/2$, then necessarily $\xi (G)\geq k+1$ holds;
moreover, $\xi (G)\geq k+2$ is valid, whenever $\left| V(G)\right| +k-1$ is
an even number. For $k=1$, we obtain a strengthening of a result of Hammer,
Hansen and Simeone, \cite{hamhansim}, which claims that $\xi (G)\geq 1$,
whenever $\alpha (G)>\left| V(G)\right| /2$. The fact that $\alpha
(G)>\left| V(G)\right| /2$ together with $\left| isol(G)\right| \neq 1$
implies $\xi (G)\geq 2$ were first established in \cite{levm1} and \cite
{levm2} for bipartite graphs and K\"{o}nig-Egerv\'{a}ry graphs respectively.
From the historical perspective it is also worth mentioning that $\xi
(T)\neq 1$ holds for any tree $T$ of order at least two, as Gunther et al., 
\cite{gun}, and independently, Zito, \cite{Zito}, have shown. Moreover, $\xi
(G)\neq 1$ is true for an arbitrary bipartite graph $G$, \cite{levm1}.

We also thoroughly investigate the special cases of K\"{o}nig-Egerv\'{a}ry
graphs and bipartite graphs, for which some sufficient conditions implying $%
\alpha (G)>\left| V(G)\right| /2$ are found. We infer that Hall's marriage
Theorem is true for K\"{o}nig-Egerv\'{a}ry graphs as well as for bipartite
graphs, and obtain a new characterization of K\"{o}nig-Egerv\'{a}ry graphs
having a perfect matching in terms of properties of $core(G)$.

\section{The main result: $\alpha (G)>(n+k-1)/2\Rightarrow $ $\xi (G)\geq
k+1 $}

Recall the following characterization of quasi-regularizable graphs, due to
Berge.

\begin{theorem}
\cite{berge2}\label{th11} A graph $G$ is quasi-regularizable if and only if $%
\left| S\right| \leq \left| N(S)\right| $ holds for any stable set $S$ of $G$%
.
\end{theorem}

\begin{corollary}
\label{cor1}If $G$ is a quasi-regularizable graph, then $\alpha (G)\leq
\left| V(G)\right| /2$.
\end{corollary}

\begin{lemma}
\label{lem1}For a graph $G$ let $S_{0}$ be a set of vertices such that 
\[
\left| S_{0}\right| =\min \{\left| S\right| :S\ is\ stable\ in\ G\ and\
\left| S\right| >\left| N(S)\right| \}. 
\]

If $isol(G)=\emptyset $, then the following assertions hold:

($\mathit{i}$) $G$ is quasi-regularizable if and only if $\left|
S_{0}\right| =0$;

($\mathit{ii}$) $G$ is not quasi-regularizable if and only if $\left|
S_{0}\right| \geq 2$.
\end{lemma}

\setlength {\parindent}{0.0cm}\textbf{Proof.} ($\mathit{i}$) It is clear,
according to Theorem \ref{th11}. Consequently, $G$ is not
quasi-regularizable if and only if $\left| S_{0}\right| \geq 1$, in fact, if
and only if $\left| S_{0}\right| \geq 2$, since $S_{0}$ is contained in some
connected component of $G$. \rule{2mm}{2mm}\setlength {\parindent}{3.45ex}

\begin{figure}[h]
\setlength{\unitlength}{1.0cm} 
\begin{picture}(5,1.6)\thicklines

  \multiput(3.5,0.5)(1,0){3}{\circle*{0.29}}
  \multiput(4.5,1.5)(1,0){2}{\circle*{0.29}}

  \put(3.5,0.5){\line(1,0){2}}
  \put(4.5,0.5){\line(0,1){1}}
  \put(4.5,0.5){\line(1,1){1}}
  \put(5.5,0.5){\line(0,1){1}}

  \put(3.5,0.84){\makebox(0,0){$a$}}
  \put(4.18,1.5){\makebox(0,0){$b$}}

  \put(4.5,0){\makebox(0,0){$G_{1}$}} 
  
  \multiput(7,0.5)(1,0){4}{\circle*{0.29}}
  \put(7,0.5){\line(1,0){3}}
  
  \multiput(8,1.5)(1,0){3}{\circle*{0.29}}
  \put(8,0.5){\line(0,1){1}}
  \put(10,0.5){\line(0,1){1}}
  
  \put(10,0.5){\line(-1,1){1}}
  \put(9,1.5){\line(1,0){1}}

  \put(9,0.84){\makebox(0,0){$c$}}
  \put(7,0.84){\makebox(0,0){$a$}}
  \put(7.67,1.5){\makebox(0,0){$b$}}
   
  \put(8.5,0){\makebox(0,0){$G_{2}$}}

\end{picture}
\caption{Non-Quasi-Regularizable graphs.}
\label{fig2}
\end{figure}
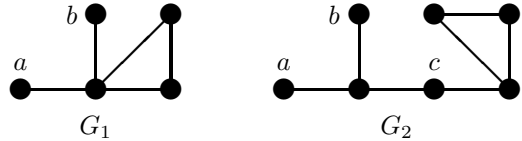

If $isol(G)\neq \emptyset $, then the above set $S_{0}$ has cardinality $%
\left| S_{0}\right| =1$. It is worth mentioning that there is a close
relation between $S_{0}$ and $core(G)$. In order to see this we need the
following characterization of a maximum stable set of a graph, due to Berge.

\begin{theorem}
\label{th2}\cite{berge2} A stable set $S$ belongs to $\Omega (G)$ if and
only if every stable set $A$ of $G$, disjoint from $S$, can be matched into $%
S$.
\end{theorem}

\begin{lemma}
\label{lem2}If $G$ is not quasi-regularizable and $S_{0}$ is a set of
vertices such that 
\[
\left| S_{0}\right| =\min \{\left| S\right| :S\ is\ stable\ in\ G\ and\
\left| S\right| >\left| N(S)\right| \}, 
\]
then $S_{0}\subseteq core(G)$.
\end{lemma}

\setlength {\parindent}{0.0cm}\textbf{Proof.} Clearly, the assertion is true
for the case that $isol(G)\neq \emptyset $, since all the isolated vertices
of $G$ are contained in any $S\in \Omega (G)$, and $S_{0}$ consists of a
single isolated vertex.\setlength {\parindent}{3.45ex}

Now, let $G$ be without isolated vertices, and suppose, on the contrary,
that there exists some $S\in \Omega (G)$, which does not include $S_{0}$.

\textit{Case 1}. $S_{0}\cap S=\emptyset $. According to Theorem \ref{th2}, $%
S_{0}$ can be matched into $S$, and this yields $\left| S_{0}\right| \leq
\left| N(S_{0})\right| $, in contradiction with the definition of $S_{0}$.

\textit{Case 2}. $S_{0}\cap S=S_{1}\neq \emptyset $. Since $S_{0}-S_{1}$ is
stable and disjoint from $S\in \Omega (G)$, Theorem \ref{th2} implies that $%
S_{0}-S_{1}$ can be matched into $S$. Hence, $\left| S_{0}-S_{1}\right| \leq
\left| N(S_{0}-S_{1})\right| $. Consequently, by $\left| S_{0}\right|
>\left| N(S_{0})\right| $ we get that $\left| S_{1}\right| >\left|
N(S_{1})\right| $, thus contradicting the minimality of $S_{0}$.

Thus, we may conclude that $S_{0}\subseteq S$, for any $S\in \Omega (G)$,
and this ensures that $S_{0}\subseteq core(G)$. \rule{2mm}{2mm}

\begin{remark}
If $G$ is not quasi-regularizable and satisfies $isol(G)=\emptyset $, then
the inclusion in Lemma \ref{lem2} may be sometimes strict. For instance, the
graph $G_{1}$ in Figure \ref{fig2} has $S_{0}=core(G_{1})=\{a,b\}$, while
for the graph $G_{2}$ in the same Figure we note that $S_{0}=\{a,b\}\subset
core(G_{2})=\{a,b,c\}$.
\end{remark}

\begin{proposition}
\label{prop9}Any $\alpha ^{+}$-stable graph free of isolated vertices is
quasi-regularizable.
\end{proposition}

\setlength {\parindent}{0.0cm}\textbf{Proof.} Assume that $G$ is a
non-quasi-regularizable $\alpha ^{+}$-stable graph with $isol(G)=\emptyset $%
. According to lemmas \ref{lem1} and \ref{lem2}, there exists some stable
set $S_{0}$ in $G$, with $\left| S_{0}\right| \geq 2$, such that $%
S_{0}\subseteq core(G)$. Consequently, we obtain that $\xi (G)\geq 2$, a
contradiction, since by virtue of Theorem \ref{th1}, $G$ must satisfy the
condition $\xi (G)\leq 1$. \rule{2mm}{2mm}\setlength {\parindent}{3.45ex}%
\newline

The restriction ''\textit{free of isolated vertices}'' in the proposition
above is essential, since no graph $G$ with $isol(G)\neq \emptyset $ can be
quasi-regularizable, but there exist $\alpha ^{+}$-stable graphs having
isolated vertices; e.g., any graph consisting of one isolated vertex and a $%
K_{n},n\geq 2$, is $\alpha ^{+}$-stable. Nevertheless, $\left|
isol(G)\right| \leq 1$ holds for any $\alpha ^{+}$-stable graph $G$.

Combining Proposition \ref{prop9} and Corollary \ref{cor1} we obtain the
following:

\begin{corollary}
\label{cor2}If $G$ is an $\alpha ^{+}$-stable graph with $isol(G)=\emptyset $%
, then $\alpha (G)\leq \left| V\left( G\right) \right| /2$.
\end{corollary}

\begin{remark}
The inequality in Corollary \ref{cor2} is not true for any $\alpha ^{+}$%
-stable graph. For instance, the graph $G$ consisting of two connected
components, namely one isolated vertex and a $C_{4}$, is $\alpha ^{+}$%
-stable and $\alpha (G)=3>\left| V\left( G\right) \right| /2=5\ /\ 2$. If $G$
consists of one isolated vertex and a $K_{n}$, $n\geq 2$, then $G$ is $%
\alpha ^{+}$-stable and $\alpha (G)=2>\left| V\left( G\right) \right| /2$,
for $n=2$, $\alpha (G)=2=\left| V\left( G\right) \right| /2$, for $n=3$,
while for $n\geq 4$, $\alpha (G)=2<\left| V\left( G\right) \right| /2$.
\end{remark}

\begin{proposition}
\label{prop4}Let $G$ be a graph with $\alpha (G)>\left| V(G)\right| /2$.
Then the following statements are equivalent:

($\mathit{i}$) $G$ is $\alpha ^{+}$-stable;

($\mathit{ii}$) $G$ is $\alpha _{1}^{+}$-stable;

($\mathit{iii}$) $G$ has a unique isolated vertex $v$, $\xi (G-v)=0$ and $%
\alpha (G-v)=(\left| V(G)\right| -1)/2$.
\end{proposition}

\setlength {\parindent}{0.0cm}\textbf{Proof.} If $\left| V(G)\right| =1$,
the result is obvious. Suppose now that $n=\left| V(G)\right| >1$.

($\mathit{i}$) $\Rightarrow $ ($\mathit{ii}$), ($\mathit{iii}$) As an $%
\alpha ^{+}$-stable graph, $G$ may have at most one isolated vertex, and
using now Corollary \ref{cor2}, we get that $G$ has exactly one isolated
vertex, say $v$. It follows that $G$ is $\alpha _{1}^{+}$-stable, $%
core(G)=\{v\}$ and $G-v$ is $\alpha _{0}^{+}$-stable with $\alpha
(G-v)=\alpha (G)-1$. Since $G-v$ is $\alpha ^{+}$-stable and has no isolated
vertices, Corollary \ref{cor2} ensures that $\alpha (G-v)\leq (n-1)/2$.
Hence, we obtain:

\[
n/2-1<\alpha (G)-1=\alpha (G-v)\leq (n-1)/2, 
\]

which implies that $n$ must be odd and $\alpha (G-v)=(n-1)/2$.%
\setlength
{\parindent}{3.45ex}

The implications ($\mathit{iii}$) $\Rightarrow $ ($\mathit{ii}$) and ($%
\mathit{ii}$) $\Rightarrow $ ($\mathit{i}$) are clear. \rule{2mm}{2mm}

\begin{corollary}
There is no graph $G$ satisfying $\alpha (G)>\left| V\left( G\right) \right|
/2$, such that either ($\mathit{i}$) $G$ is $\alpha _{1}^{+}$-stable and of
even order, or ($\mathit{ii}$) $G$ is $\alpha _{0}^{+}$-stable.
\end{corollary}

\begin{corollary}
\label{cor5}If $\alpha (G)>\left| V\left( G\right) \right| /2$ and $\xi
(G)=1 $, then $\left| V(G)\right| \equiv 1$ \textit{mod} $2$.
\end{corollary}

\begin{theorem}
\label{th3}Let $G$ be a graph with $\alpha (G)>\left| V\left( G\right)
\right| /2$. Then the following statements are true:

($\mathit{i}$) if $\left| isol(G)\right| \neq 1$, then $\xi (G)\geq 2$;

($\mathit{ii}$) if $\xi (G)=1$, then $\left| isol(G)\right| =1$;

($\mathit{iii}$) $\xi (G)=0$ never holds.
\end{theorem}

\setlength {\parindent}{0.0cm}\textbf{Proof.} ($\mathit{i}$) The result is
clear, whenever $\left| isol(G)\right| \geq 2$, because $\xi (G)\geq \left|
isol(G)\right| $. Let $isol(G)=\emptyset $, and suppose, on the contrary,
that $\xi (G)\leq 1$. By Theorem \ref{th1}, we infer that $G$ is an $\alpha
^{+}$-stable graph. Hence, Corollary \ref{cor2} ensures that $\alpha (G)\leq
\left| V\left( G\right) \right| /2$, in contradiction with the premise on $G$%
. Therefore, $\xi (G)\geq 2$ is true.\setlength
{\parindent}{3.45ex}

($\mathit{ii}$) If $\xi (G)=1$, then $\left| isol(G)\right| =1$, since
otherwise, according to ($\mathit{i}$), $G$ must satisfy $\xi (G)\geq 2$.

($\mathit{iii}$) If $\xi (G)=0$, then $G$ is $\alpha ^{+}$-stable. Since $%
\alpha (G)>\left| V\left( G\right) \right| /2$, Proposition \ref{prop4}
implies that $G$ must be $\alpha _{1}^{+}$-stable, i.e., $\xi (G)=1$, in
contradiction with the assumption on $G$. \rule{2mm}{2mm}

\begin{remark}
The condition $\xi (G)\geq 2$ is not strong enough to ensure $\alpha (G)>n/2$%
. Moreover, for any positive integer $k$, one can choose an arbitrarily
large positive integer $p$ and build a graph $G=(V,E)$, with $n=k+p$
vertices, such that $\xi (G)=k$ and $\alpha (G)<n/2$. For instance, such a
graph $G=(V,E)$ is defined by: 
\[
V=\{x_{i}:1\leq i\leq k\}\cup \{y_{i}:1\leq i\leq p\},\ E=\{x_{i}y_{1}:1\leq
i\leq k\}\cup \{y_{i}y_{j}:1\leq i<j\leq p\}. 
\]
\end{remark}

Theorem \ref{th3} is a strengthening of the following result, due to Hammer,
Hansen and Simeone:

\begin{corollary}
\cite{hamhansim}\label{cor4} If $G$ has $\alpha (G)>\left| V\left( G\right)
\right| /2$, then $\xi (G)\geq 1$.
\end{corollary}

\begin{proposition}
\label{prop6}If $G=(V,E)$ and $H=G[V-N[core(G)]]$, then the following
assertions are true:

($\mathit{i}$) $H$ has no isolated vertices;

($\mathit{ii}$) $\alpha (H)=\alpha (G)-\xi (G)$;

($\mathit{iii}$) $S_{H}\in \Omega (H)$ if and only if there is $S_{G}\in
\Omega (G)$, such that $S_{H}=S_{G}\cap V(H)$;

($\mathit{iv}$) $H$ is $\alpha _{0}^{+}$-stable;

($\mathit{v}$) $\left| S-core(G)\right| \leq \left| N(S)-N(core(G))\right| $
holds for any $S\in \Omega (G)$.
\end{proposition}

\setlength {\parindent}{0.0cm}\textbf{Proof.} If $core(G)=\emptyset $, then $%
H=G$ and all the assertions are valid.\setlength {\parindent}{3.45ex} Assume
that $core(G)\neq \emptyset $. Let $S\in \Omega (G),A=S-core(G)$ and $%
B=V-S-N(core(G))$.

($\mathit{i}$) Suppose, on the contrary, that $H$ has an isolated vertex $v$.

\textit{Case 1}. $v\in B$. Then $(\{v\},core(G))=(\{v\},A)=\emptyset $,
because $v\notin N(core(G))$ and it is isolated in $H$. Hence, $S\cup \{v\}$
is stable in $G$, contradicting the maximality of $S$.

\textit{Case 2}. $v\in A$. Let $W\in \Omega (G),W\neq S$. Then $%
core(G)\subset W,W\cap B\neq \emptyset ,W\cap N(core(G))=\emptyset $ and
maybe $W\cap A\neq \emptyset $. Hence, $v$ is adjacent to no vertex in $W$,
and since $W$ is a maximum stable set in $G$, we infer that $v\in W$. It
follows that $v\in core(G)$, in contradiction with the fact that $v\in
A=S-core(G)$.

($\mathit{ii}$) Since $A$ is stable, it follows that $\alpha (G)-\xi
(G)=\left| A\right| \leq \alpha (H)$. Suppose, on the contrary, that $\alpha
(H)=\left| S_{H}\right| >\left| A\right| $, with $S_{H}\in \Omega (H)$.
Hence, we get that: $(S_{H}\cap A,core(G))=(S_{H}\cap B,core(G))=\emptyset $%
, because $S_{H}\cap A\subset A\subset S$ and also $B\cap
N(core(G))=\emptyset $. Consequently, $S_{H}\cup core(G)$ is stable in $G$
and $\left| S_{H}\cup core(G)\right| =\left| S_{H}\right| +\xi (G)>\left|
A\right| +\xi (G)=\alpha (G)$, in contradiction with definition of $\alpha
(G)$.

($\mathit{iii}$) If $S_{H}\in \Omega (H)$, then $(S_{H}\cap
B,core(G))=\emptyset $, because $B\cap N(core(G))=\emptyset $, and evidently 
$(S_{H}\cap A,core(G))=\emptyset $. Hence, $S_{G}=core(G)\cup S_{H}$ is
stable in $G$ and since, by ($\mathit{ii}$), $\alpha (H)=\left| S_{H}\right|
=\alpha (G)-\xi (G)$, we infer that $S_{G}\in \Omega (G)$.

Conversely, if $S_{H}=S_{G}\cap V(H)$ and $S_{G}\in \Omega (G)$, then $S_{H}$
is stable in $H$ as well, and $\left| S_{H}\right| =\left| S_{G}\right| -\xi
(G)=\alpha (G)-\xi (G)=\alpha (H)$. Therefore, $S_{H}\in \Omega (H)$.

($\mathit{iv}$) According to ($\mathit{iii}$), it follows that $%
core(H)=\emptyset $, i.e., $H$ is $\alpha _{0}^{+}$-stable.

($\mathit{v}$) $H$ is $\alpha _{0}^{+}$-stable and has no isolated vertices.
Hence, Corollary \ref{cor2} implies that $\alpha (H)\leq \left| V(H)\right|
/2$, i.e., $\left| S-core(G)\right| =\left| A\right| \leq \left| B\right|
=\left| N(S)-N(core(G))\right| .$ \rule{2mm}{2mm}

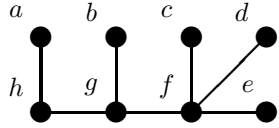
\begin{figure}[h]
\setlength{\unitlength}{1.0cm} 
\begin{picture}(5,1.1)\thicklines
  
  \multiput(5.5,0)(1,0){4}{\circle*{0.29}}
  \multiput(5.5,1)(1,0){4}{\circle*{0.29}}
  \put(5.5,0){\line(1,0){3}}
  
  \put(5.17,0.34){\makebox(0,0){$h$}}
  \put(6.17,0.34){\makebox(0,0){$g$}}
  \put(7.17,0.34){\makebox(0,0){$f$}}
  \put(8.25,0.34){\makebox(0,0){$e$}}

  \multiput(5.5,0)(1,0){3}{\line(0,1){1}}
  \put(7.5,0){\line(1,1){1}} 

  \put(5.17,1.34){\makebox(0,0){$a$}}
  \put(6.17,1.34){\makebox(0,0){$b$}}
  \put(7.17,1.34){\makebox(0,0){$c$}}
  \put(8.17,1.34){\makebox(0,0){$d$}}

\end{picture}
\caption{{A counterexample to $\left| S-A\right| \leq \left|
N(S)-N(A)\right| $.}}
\label{fig34}
\end{figure}

\begin{remark}
The inequality $\left| S-A\right| \leq \left| N\left( S\right) -N\left(
A\right) \right| $ is not generally true for any subset $A$ of a maximum
stable set $S$. For instance, for the graph in Figure \ref{fig34}, one can
take $S=\{a,b,c,d,e\},A=\{b,c,d\}$, and get $\left| S-A\right| =2>\left|
N(S)-N(A)\right| =1.$
\end{remark}

By Corollary \ref{cor1}, if $\alpha (G)>\left| V\left( G\right) \right| /2$,
then $G$ is not quasi-regularizable. From the point of view of Berge's
Theorem \ref{th11}, an elementary obstacle to being quasi-regularizable is
the presence of a stable set $S$ in $G$ satisfying $\left| S\right| >\left|
N(S)\right| $. An interesting question is how to present at least one such
stable set for the graph $G$. The following result gives a canonical
procedure for constructing such a stable set, when $\alpha (G)>\left|
V\left( G\right) \right| /2$.

\begin{theorem}
\label{th4}If a graph $G$ has $\alpha (G)>\left| V\left( G\right) \right| /2$%
, then $\left| core(G)\right| >\left| N(core(G))\right| $.
\end{theorem}

\setlength {\parindent}{0.0cm}\textbf{Proof.} Let $S\in \Omega (G),$ $%
A=S-core(G),B=V-S-N(core(G))$, and also $H=G[V-N[core(G)])]$. By Proposition 
\ref{prop6}, $\alpha (H)=\alpha (G)-\xi (G)=\left| A\right| $. Suppose, on
the contrary, that $\xi (G)\leq \left| N(core(G))\right| $. Since $\alpha
(G)>\left| V\left( G\right) \right| /2$ we get $\alpha (G)=\xi (G)+\left|
A\right| >\left| N(core(G))\right| +\left| B\right| $, which together with $%
\xi (G)\leq \left| N(core(G))\right| $ gives $\left| A\right| =\alpha
(G)-\xi (G)>\left| V\left( G\right) \right| /2-\left| N(core(G))\right|
=\left| B\right| $. Therefore, it follows $\alpha (H)=\left| A\right|
>\left| V(H)\right| /2$. Using again Proposition \ref{prop6}, we infer that $%
H$ is $\alpha _{0}^{+}$-stable and $isol(H)=\emptyset $. According to
Corollary \ref{cor2}, it follows that $\alpha (H)\leq \left| V(H)\right| /2$%
, which contradicts the premise of the theorem. \rule{2mm}{2mm}%
\setlength
{\parindent}{3.45ex}\newline

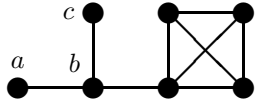
\begin{figure}[h]
\setlength{\unitlength}{1.0cm} 
\begin{picture}(5,1)\thicklines
  
  \multiput(5.5,0)(1,0){4}{\circle*{0.29}}
  \multiput(6.5,1)(1,0){3}{\circle*{0.29}}
  \put(5.5,0){\line(1,0){3}}
   \put(7.5,1){\line(1,0){1}}
  
  \put(5.5,0.34){\makebox(0,0){$a$}}
  \put(6.25,0.34){\makebox(0,0){$b$}}

  \multiput(6.5,0)(1,0){3}{\line(0,1){1}}
  \put(7.5,0){\line(1,1){1}}
   \put(7.5,1){\line(1,-1){1}} 

  \put(6.17,1){\makebox(0,0){$c$}}

\end{picture}
\caption{{$\xi (G)>\left| N(core(G))\right| $ and $\alpha (G)<\left| V\left(
G\right) \right| /2.$}}
\label{fig4}
\end{figure}

The converse of Theorem \ref{th4} is not generally true. The graph in Figure 
\ref{fig4} illustrates this assertion.

\begin{proposition}
\label{prop7}If a graph $G$ has $\alpha (G)>\left| V\left( G\right) \right|
/2$, and $\xi (G)\leq k$, for some $k\geq 2$, then $\alpha (G)\leq (\left|
V(G)\right| +k-1)/2$.
\end{proposition}

\setlength {\parindent}{0.0cm}\textbf{Proof.} The result is clear for $%
\left| V(G)\right| =1$.

If $\left| V(G)\right| >1$, let us denote $H=G[V(G)-N[core(G)]]$. Then we
obtain: 
\[
\left| V(H)\right| =\left| V(G)\right| -\left| core(G)\right| -\left|
N(core(G))\right| =\left| V(G)\right| -\xi (G)-\left| N(core(G))\right| . 
\]

By Proposition \ref{prop6}, $H$ is $\alpha _{0}^{+}$-stable and $%
isol(H)=\emptyset $. According to Corollary \ref{cor2}, it follows that $%
\alpha (H)\leq \left| V(H)\right| /2\leq (\left| V(G)\right| -\xi (G)-1)/2$,
and consequently, we infer that: $\alpha (G)=\alpha (H)+\xi (G)\leq (\left|
V(G)\right| -\xi (G)-1)/2+\xi (G)=(\left| V(G)\right| +\xi (G)-1)/2\leq
(\left| V(G)\right| +k-1)/2$. \rule{2mm}{2mm}\setlength {\parindent}{3.45ex}

\begin{theorem}
\label{th5}If $\left| isol(G)\right| \neq 1$ and $\alpha (G)>(\left|
V(G)\right| +k-1)/2$, for some $k\geq 1$, then the following assertions are
true:

($\mathit{i}$) $\xi (G)\geq k+1$, whenever $\left| V(G)\right| +k-1$ is odd;

($\mathit{ii}$) $\xi (G)\geq k+2$, whenever $\left| V(G)\right| +k-1$ is
even.
\end{theorem}

\setlength {\parindent}{0.0cm}\textbf{Proof.} Clearly, $n=\left| V(G)\right|
\neq 1$. First we show that ($\mathit{i}$) is always valid. For $k=1$ this
assertion is exactly the claim of Theorem \ref{th3}($\mathit{i}$). Suppose
that $k\geq 2$, and $\alpha (G)>(n+k-1)/2\geq $ $n/2$. If $\xi (G)\leq k$,
then by Proposition \ref{prop7}, it follows that $\alpha (G)\leq (n+k-1)/2$,
in contradiction with the premises on $\alpha (G)$. Therefore, $\xi (G)\geq
k+1$ is true for any $k\geq 1$. Assume that ($\mathit{ii}$) is not valid.
Therefore we get that $\xi (G)\leq k+1$, and because $\alpha
(G)>(n+k-1)/2\geq $ $n/2$, Proposition \ref{prop7} implies that $\alpha
(G)\leq (n+k+1-1)/2=(n+k)/2$. Hence we obtain the following contradiction:

\[
(n+k-1)/2<\alpha (G)\leq (n+k)/2, 
\]

since $\alpha (G)$ must be a positive integer. \rule{2mm}{2mm}%
\setlength
{\parindent}{3.45ex}\newline

For $k=1$, we obtain the following strengthening of Corollary \ref{cor4} due
to Hammer, Hansen and Simeone, and of Theorem \ref{th3}($\mathit{i}$):

\begin{corollary}
If $\left| isol(G)\right| \neq 1$ and $\alpha (G)>\left| V(G)\right| /2$,
then $\xi (G)\geq 2$, whenever $\left| V(G)\right| $ is odd, and $\xi
(G)\geq 3$, if $\left| V(G)\right| $ is even.
\end{corollary}

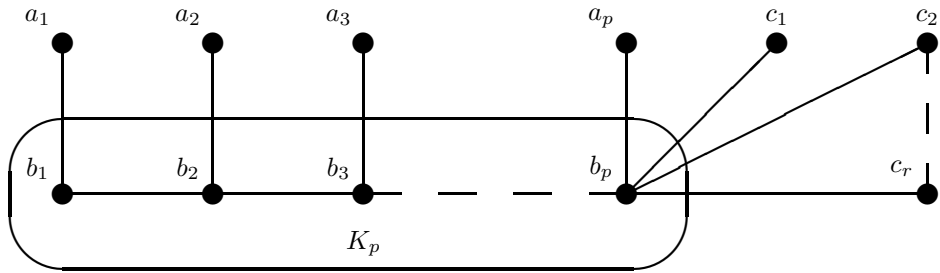
\begin{figure}[h]
\setlength{\unitlength}{1.0cm} 
\begin{picture}(5,3.4)\thicklines

  \put(4.8,1){\oval(9.0,2)}

  \multiput(1,1)(2,0){3}{\circle*{0.29}}
  \multiput(1,3)(2,0){3}{\circle*{0.29}}

  \put(1,1){\line(1,0){4}}

  \multiput(1,1)(2,0){3}{\line(0,1){2}}

  \put(0.67,1.34){\makebox(0,0){$b_{1}$}}
  \put(2.67,1.34){\makebox(0,0){$b_{2}$}}
  \put(4.67,1.34){\makebox(0,0){$b_{3}$}}
  \put(8.16,1.34){\makebox(0,0){$b_{p}$}}

  \multiput(5,1)(1,0){4}{\line(1,0){0.5}}

  \multiput(8.5,1)(4,0){2}{\circle*{0.29}}

  \put(12.16,1.34){\makebox(0,0){$c_{r}$}}

  \multiput(8.5,3)(2,0){3}{\circle*{0.29}}

  \put(8.5,1){\line(1,0){4}}
  \put(8.5,1){\line(0,1){2}}

  \put(8.5,1){\line(1,1){2}}
  \put(8.5,1){\line(2,1){4}}

  \multiput(12.5,1)(0,0.8){3}{\line(0,1){0.4}}

  \put(0.67,3.34){\makebox(0,0){$a_{1}$}}
  \put(2.67,3.34){\makebox(0,0){$a_{2}$}}
  \put(4.67,3.34){\makebox(0,0){$a_{3}$}}
  \put(8.16,3.34){\makebox(0,0){$a_{p}$}}

  \put(10.5,3.34){\makebox(0,0){$c_{1}$}}
  \put(12.5,3.34){\makebox(0,0){$c_{2}$}}

  \put(5,0.34){\makebox(0,0){$K_{p}$}}
 
\end{picture}
\caption{{$\alpha (G)=p+r,\ \xi (G)=r+1$.}}
\label{fig45}
\end{figure}

\begin{remark}
The graph $G$ in Figure \ref{fig45}\emph{\ }shows that the bounds in Theorem 
\ref{th5} are tight. We may consider either $k=r-1$ or $k=r$, and
correspondingly, $n+k-1$ will be even or odd.
\end{remark}

\section{K\"{o}nig-Egerv\'{a}ry graphs and bipartite graphs}

In the sequel we deal with K\"{o}nig-Egerv\'{a}ry graphs, for which we show
that the converse of Theorem \ref{th4} is also true. It is known that $%
\lfloor n/2\rfloor +1\leq \alpha (G)+\mu (G)\leq n$ holds for any graph $G$
with $n$ vertices.

\begin{lemma}
\label{lem4}If $G$ is a K\"{o}nig-Egerv\'{a}ry graph, then $\alpha (G)\geq
\mu (G)$.
\end{lemma}

\setlength {\parindent}{0.0cm}\textbf{Proof.} Since $\mu (G)\leq \left|
V(G)\right| /2$ holds for any graph $G$, and $\left| V(G)\right| =\alpha
(G)+\mu (G)$ is true for $G$ a K\"{o}nig-Egerv\'{a}ry graph, then clearly
follows that $\alpha (G)\geq \mu (G)$ is valid in our premise on $G$. \rule%
{2mm}{2mm}\setlength {\parindent}{3.45ex}

\begin{corollary}
\label{cor3}If $G$ is a K\"{o}nig-Egerv\'{a}ry graph, then the following
statements are true:

($\mathit{i}$) $\alpha (G)\geq \left| V(G)\right| /2$;

($\mathit{ii}$) $\alpha (G)=\left| V(G)\right| /2$ if and only if $G$ has a
perfect matching;

($\mathit{iii}$) $\alpha (G)>\left| V(G)\right| /2$, whenever $isol(G)\neq
\emptyset $.
\end{corollary}

Combining Proposition \ref{prop4} and Corollary \ref{cor3} we may conclude
that:

\begin{proposition}
\label{prop3}If $G$ is an $\alpha ^{+}$-stable K\"{o}nig-Egerv\'{a}ry graph
with $\left| isol\left( G\right) \right| \neq 1$, i.e., $isol\left( G\right)
\neq \emptyset $, then $\alpha (G)=\left| V(G)\right| /2$.
\end{proposition}

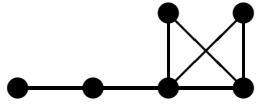
\begin{figure}[h]
\setlength{\unitlength}{1.0cm} 
\begin{picture}(5,1)\thicklines
  
  \multiput(5.5,0)(1,0){4}{\circle*{0.29}}
  \multiput(7.5,1)(1,0){2}{\circle*{0.29}}
  \put(5.5,0){\line(1,0){3}}
  
  \multiput(7.5,0)(1,0){2}{\line(0,1){1}}
  \put(7.5,0){\line(1,1){1}}
   \put(7.5,1){\line(1,-1){1}} 

\end{picture}
\caption{{A Koenig-Egervary graph $G$ with $\alpha (G)=\left| V\left(
G\right) \right| /2$ and $\xi (G)\geq 2$.}}
\label{fig3}
\end{figure}

\begin{remark}
There are non-$\alpha ^{+}$-stable K\"{o}nig-Egerv\'{a}ry graphs $G$,
without isolated vertices, such that $\alpha (G)=\left| V\left( G\right)
\right| /2$. For instance, the graph $G$ in Figure \ref{fig3} has $\alpha
(G)=3=\left| V\left( G\right) \right| /2$ and $\xi (G)=2$, i.e., $G$ is not $%
\alpha ^{+}$-stable. Notice that $G$ has perfect matchings.
\end{remark}

\begin{proposition}
\label{prop5}If $G$ is a K\"{o}nig-Egerv\'{a}ry graph, then the following
statements are equivalent:

($\mathit{i}$) $\alpha (G)>\left| V\left( G\right) \right| /2$;

($\mathit{ii}$) $G$ has no perfect matching;

($\mathit{iii}$) $G$ is non-quasi-regularizable;

($\mathit{iv}$) $\xi (G)=\left| core(G)\right| >\left| N(core(G))\right| $.
\end{proposition}

\setlength {\parindent}{0.0cm}\textbf{Proof.} Theorem \ref{th4} ensures that
($\mathit{i}$) $\Rightarrow $($\mathit{iv}$); the step ($\mathit{iv}$) $%
\Rightarrow $($\mathit{iii}$) is clear by Theorem \ref{th11}; the
implication ($\mathit{iii}$) $\Rightarrow $($\mathit{ii}$) follows from
Theorem \ref{th11}; and ($\mathit{ii}$) $\Rightarrow $ ($\mathit{i}$) is
true according to Corollary \ref{cor3}. \rule{2mm}{2mm}%
\setlength
{\parindent}{3.45ex}

\begin{remark}
In general, ($\mathit{ii}$) implies neither ($\mathit{i}$) nor ($\mathit{iii}
$); e.g., the graph $K_{3}$. The graph in Figure \ref{fig4} shows that ($%
\mathit{i}$) does not always follow from ($\mathit{iii}$) or from ($\mathit{%
iv}$).
\end{remark}

\begin{remark}
If $G$ is a K\"{o}nig-Egerv\'{a}ry graph with $isol(G)\neq \emptyset $, then
all the assertions in Proposition \ref{prop5} are valid, since ($\mathit{i}$%
) follows from Corollary \ref{cor3}.
\end{remark}

\begin{theorem}
\label{th6}If $G$ is a K\"{o}nig-Egerv\'{a}ry graph, then the following
statements are equivalent:

($\mathit{i}$) $\alpha (G)=\left| V\left( G\right) \right| /2$;

($\mathit{ii}$) $G$ has a perfect matching;

($\mathit{iii}$) $G$ is quasi-regularizable;

($\mathit{iv}$) $\left| core(G)\right| \leq \left| N(core(G))\right| $.
\end{theorem}

\setlength {\parindent}{0.0cm}\textbf{Proof.} It follows from Proposition 
\ref{prop5} and Lemma \ref{lem4}. \rule{2mm}{2mm}%
\setlength
{\parindent}{3.45ex}

\begin{remark}
There are quasi-regularizable graphs without perfect matching; e.g., the
graph in Figure \ref{fig56}.
\end{remark}

\begin{figure}[h]
\setlength{\unitlength}{1.0cm} 
\begin{picture}(5,2)\thicklines
  
  \multiput(3,0)(2,0){2}{\circle*{0.29}}
  \multiput(9,0)(2,0){2}{\circle*{0.29}}

  \put(3,0){\line(1,0){8}}

  \multiput(5,1)(4,0){2}{\circle*{0.29}}

  \put(7,2){\circle*{0.29}}
  
  \multiput(5,0)(4,0){2}{\line(0,1){1}}
  \put(5,1){\line(1,0){4}} 
  \put(5,1){\line(2,1){2}} 
  \put(5,0){\line(4,1){4}}
  \put(5,0){\line(1,1){2}}  
  \put(9,0){\line(-4,1){4}}
  \put(9,0){\line(-1,1){2}}
  \put(9,1){\line(-2,1){2}}

\end{picture}
\caption{{$G$ is quasi-regularizable, but has no perfect matching. }}
\label{fig56}
\end{figure}

\begin{remark}
The equivalence ($\mathit{ii}$) $\Leftrightarrow $ ($\mathit{iii}$) in
Theorem \ref{th6} (and comparison with Theorem \ref{th11}) shows that Hall's
marriage Theorem is also true for K\"{o}nig-Egerv\'{a}ry graphs. Moreover,
according to ($\mathit{iv}$), it is sufficient to check Hall's condition
only for one specific stable set, namely for $core(G)$.
\end{remark}

\begin{remark}
The graph in Figure \ref{fig57} shows that there exist
non-quasi-regularizable graphs satisfying the inequality $\left|
core(G)\right| \leq \left| N(core(G))\right| $.
\end{remark}

\begin{figure}[h]
\setlength{\unitlength}{1.0cm} 
\begin{picture}(5,2)\thicklines

  \multiput(3,0)(2,0){2}{\circle*{0.29}}
  \multiput(3,2)(2,0){2}{\circle*{0.29}} 
  \multiput(7,0)(4,0){2}{\circle*{0.29}}
  \multiput(7,1)(4,0){2}{\circle*{0.29}}

  \put(9,2){\circle*{0.29}}
  \put(3,0){\line(1,0){8}}

  \put(3,0){\line(0,1){2}}
  \put(3,0){\line(1,1){2}}  
  
  \put(5,0){\line(2,1){2}} 
  \put(5,0){\line(6,1){6}}

  \multiput(7,0)(4,0){2}{\line(0,1){1}}
  
  \put(7,1){\line(1,0){4}} 
  \put(7,1){\line(2,1){2}} 
  \put(7,0){\line(4,1){4}}
  \put(7,0){\line(1,1){2}}  
 
  \put(11,0){\line(-4,1){4}}
  \put(11,0){\line(-1,1){2}}
  \put(11,1){\line(-2,1){2}}

\end{picture}
\caption{{$G$ is not quasi-regularizable, and $\left| core\left( G\right)
\right| =3\leq \left| N(core(G))\right| =4.$}}
\label{fig57}
\end{figure}

\begin{proposition}
\label{prop8}If $G=\left( A,B,E\right) $ is a bipartite graph with $\left|
isol\left( G\right) \right| \neq 1$, then either ($\mathit{i}$) $\xi (G)\geq
2$, or ($\mathit{ii}$) $\xi (G)=0$ and $A,B\in \Omega (G)$.
\end{proposition}

\setlength {\parindent}{0.0cm}\textbf{Proof.} By Corollary \ref{cor3}, $G$
satisfies $\alpha (G)\geq \left| V\left( G\right) \right| /2$. If $\alpha
(G)>\left| V\left( G\right) \right| /2$, then $\left| isol\left( G\right)
\right| \neq 1$ assures that $\xi (G)\geq 2$ holds, by Theorem \ref{th3}. If 
$\alpha (G)=\left| V\left( G\right) \right| /2$, then Corollary \ref{cor3}($%
\mathit{ii}$) implies that $\left| A\right| =\left| B\right| =\left| V\left(
G\right) \right| /2$. Hence, $A,B\in \Omega (G)$ and evidently $\xi (G)=0$. 
\rule{2mm}{2mm}\setlength {\parindent}{3.45ex}

\begin{figure}[h]
\setlength{\unitlength}{1.0cm} 
\begin{picture}(5,2.5)\thicklines

  \multiput(1,0)(2,0){3}{\circle*{0.29}}
  \multiput(1,2)(2,0){3}{\circle*{0.29}}

  \put(1,0){\line(1,0){4}}

  \multiput(1,0)(2,0){3}{\line(0,1){2}}

  \put(0.67,0.34){\makebox(0,0){$b_{1}$}}
  \put(2.67,0.34){\makebox(0,0){$b_{2}$}}
  \put(4.67,0.34){\makebox(0,0){$b_{3}$}}
  \put(8.16,0.34){\makebox(0,0){$b_{p}$}}

  \multiput(5,0)(1,0){4}{\line(1,0){0.5}}

  \multiput(8.5,0)(2,0){3}{\circle*{0.29}}

  \put(10.5,0.34){\makebox(0,0){$c_{1}$}}
  \put(12.5,0.34){\makebox(0,0){$c_{2}$}}

  \put(8.5,2){\circle*{0.29}}

  \put(8.5,0){\line(1,0){4}}
  \put(8.5,0){\line(0,1){2}}
  \put(10.5,0){\line(-1,1){2}}

  \put(0.67,2.34){\makebox(0,0){$a_{1}$}}
  \put(2.67,2.34){\makebox(0,0){$a_{2}$}}
  \put(4.67,2.34){\makebox(0,0){$a_{3}$}}
  \put(8.16,2.34){\makebox(0,0){$a_{p}$}}
 
\end{picture}
\caption{{A Koenig-Egervary graph $G$ of order $k=2p+2$ with $\xi (G)=1$.}}
\label{fig77}
\end{figure}
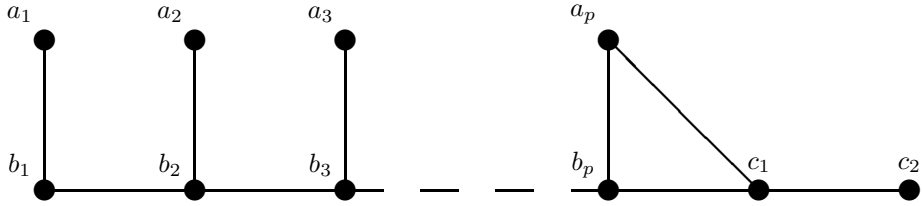

\begin{remark}
The above result is not true for all K\"{o}nig-Egerv\'{a}ry graphs.
Moreover, for any even positive integer $k$ there exists a
K\"{o}nig-Egerv\'{a}ry graph $G$ of order $k$ with $\xi (G)=1$. This
observation is illustrated in Figure \ref{fig77}, where $core(G)=\left\{
c_{2}\right\} $.

If the order of a K\"{o}nig-Egerv\'{a}ry graph $G$ is odd, then $\alpha
(G)>\left| V\left( G\right) \right| /2$ and consequently $\xi (G)\geq 2$.
However, for every odd positive integer $k$ there exists a graph $G$ of size 
$k$ with $\xi (G)=1$. For instance, the graph $H=(V(G)\cup
\{c_{3}\},E(G)\cup \{a_{1}c_{3},b_{1}c_{3}\})$ is of odd order $k=2p+3$ and $%
core(H)=\left\{ c_{2}\right\} $.
\end{remark}

\begin{theorem}
If $G$ is a bipartite graph with $\xi (G)=1$, then $\left| isol(G)\right| =1$
and $\left| V(G)\right| \equiv 1$ \textit{mod} $2$.
\end{theorem}

\setlength {\parindent}{0.0cm}\textbf{Proof.} Proposition \ref{prop8}
ensures that $\left| isol(G)\right| =1$. Then, Corollary \ref{cor3} implies
that $\alpha (G)>\left| V\left( G\right) \right| /2$, and consequently, $%
\left| V(G)\right| \equiv 1$ \textit{mod} $2$ is now true, according to
Corollary \ref{cor5}. \rule{2mm}{2mm}\setlength {\parindent}{3.45ex}\newline

As a consequence of Proposition \ref{prop8} we obtain the following
characterization of $\alpha ^{+}$-stable bipartite graphs.

\begin{corollary}
\label{cor6}\cite{levm} For a bipartite graph $G$ with $\left|
isol(G)\right| \neq 1$, the following assertions are equivalent:

($\mathit{i}$) $G$ is $\alpha ^{+}$-stable;

($\mathit{ii}$) $G$ has a perfect matching;

($\mathit{iii}$) $G$ possesses two maximum stable sets that partition its
vertex set.
\end{corollary}

\setlength {\parindent}{0.0cm}\textbf{Proof.} ($\mathit{i}$) $\Rightarrow $ (%
$\mathit{ii}$), ($\mathit{iii}$) By Proposition \ref{prop3}, we get that $%
\alpha (G)=\left| V\left( G\right) \right| /2$. Hence, Theorem \ref{th6}
implies ($\mathit{ii}$), and ($\mathit{iii}$) follows from Proposition \ref
{prop8}.\setlength {\parindent}{3.45ex}

($\mathit{iii}$) $\Rightarrow $ ($\mathit{i}$) It is evident, by Theorem \ref
{th1}.

($\mathit{ii}$) $\Rightarrow $ ($\mathit{iii}$) If $G$ has a perfect
matching, then Theorem \ref{th6} assures that $\alpha (G)=\left| V\left(
G\right) \right| /2$, and this leads (see the proof of Proposition \ref
{prop8}) to $A,B\in \Omega (G)$. \rule{2mm}{2mm}\newline

This result generalizes the corresponding statement for trees, given in \cite
{gun}.

Clearly, Proposition \ref{prop5} is true for bipartite graphs as well, since
they are K\"{o}nig-Egerv\'{a}ry graphs. Moreover, using Corollary \ref{cor6}%
, we obtain a stronger form of this assertion.

\begin{proposition}
If $G$ is a bipartite graph with $\left| isol\left( G\right) \right| \neq 1$%
, then 
\[
\alpha (G)>\left| V\left( G\right) \right| /2\ if\ and\ only\ if\ \xi
(G)\geq 2. 
\]
\end{proposition}

\setlength {\parindent}{0.0cm}\textbf{Proof.} According to Theorem \ref{th3}
($\mathit{i}$), $\left| isol\left( G\right) \right| \neq 1$ and $\alpha
(G)>\left| V\left( G\right) \right| /2$ implies that $\xi (G)\geq 2$ is true
for general graphs.\setlength {\parindent}{3.45ex}

Conversely, for $G$ a bipartite graph with $\left| isol\left( G\right)
\right| \neq 1$, we have to show that if $\xi (G)\geq 2$, then also $\alpha
(G)>\left| V\left( G\right) \right| /2$. By Theorem \ref{th1}, it follows
that $G$ is not $\alpha ^{+}$-stable, and then Corollary \ref{cor6} implies
that $G$ has no perfect matching. Hence, by Corollary \ref{cor3}, we get
that $\alpha (G)>\left| V\left( G\right) \right| /2$. \rule{2mm}{2mm}

\section{Conclusions}

In this paper we have presented a number of relations connecting $\alpha (G)$
and $\xi (G)$, whenever $\alpha (G)>\left| V(G)\right| /2$, with emphasis on
K\"{o}nig-Egerv\'{a}ry graphs and bipartite graphs, for which some
sufficient conditions are also necessary. It would be interesting to see
which of these results can be transferred to graphs satisfying $\alpha
(G)=\left| V(G)\right| /2$ or less. The special case of
K\"{o}nig-Egerv\'{a}ry graphs could offer a promising start, since for them
the condition $\alpha (G)\leq \left| V(G)\right| /2$ is equivalent to $%
\alpha (G)=\left| V(G)\right| /2$. We also have shown that a necessary and
sufficient condition for a K\"{o}nig-Egerv\'{a}ry graph $G$ to have a
perfect matching reads as $\left| core(G)\right| \leq \left|
N(core(G))\right| $. The following challenging problem seems to be of
algorithmic nature: how to find $core(G)$, at least for
K\"{o}nig-Egerv\'{a}ry graphs.

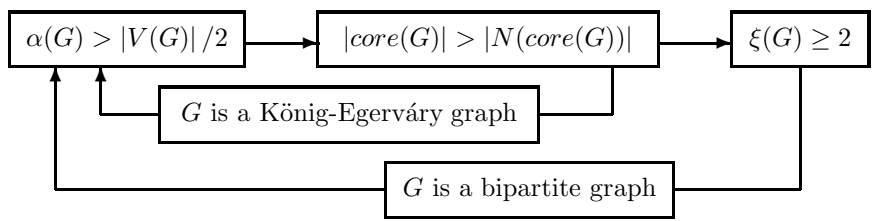
\begin{figure}[h]
\setlength{\unitlength}{1.0cm} 
\begin{picture}(5,3)\thicklines

\put(1,2){\framebox(3.1,0.7){$\alpha (G)>\left| V(G)\right| /2$}}
\put(4.1,2.3){\vector(1,0){1}}

\put(5.1,2){\framebox(4.5,0.7){$\left| core(G)\right| >\left| N(core(G))\right| $}}
\put(9.6,2.3){\vector(1,0){1}}

\put(9,2){\line(0,-1){0.65}}
\put(9,1.35){\line(-1,0){1}}

\put(3,1.35){\line(-1,0){0.8}}
\put(2.2,1.35){\vector(0,1){0.65}}

\put(10.6,2){\framebox(1.8,0.7){$\xi (G)\geq 2$}}

\put(11.5,2){\line(0,-1){1.65}}
\put(11.5,0.35){\line(-1,0){1.7}}

\put(6,0.35){\line(-1,0){4.4}}
\put(1.6,0.35){\vector(0,1){1.65}}

\put(3,1){\framebox(5,0.7){$G$ is a K\"{o}nig-Egerv\'{a}ry graph}}

\put(6,0){\framebox(3.8,0.7){$G$ is a bipartite graph}}

 \end{picture}
\caption{A scheme of interconnections between the main findings of the
paper. }
\label{fig99}
\end{figure}

Fugure \ref{fig99} suggests the following questions: what must be added to $%
\xi (G)\geq 2$, in order to get $\alpha (G)>\left| V(G)\right| /2$ for, at
least, K\"{o}nig-Egerv\'{a}ry graphs; when the inequality $\left|
core(G)\right| >\left| N(core(G))\right| $ can be obtained from $\xi (G)\geq
2$?

\section{Acknowledgment}

The authors thank Professor Uri Peled for his careful reading of this
manuscript, and also for his suggestions and comments.

\end{document}